\documentclass[12pt]{amsart}
\usepackage{amsfonts,amssymb,verbatim,amsmath,amsthm,latexsym,textcomp,amscd,mathtools,gensymb,subcaption}
\usepackage{latexsym,amsfonts,amssymb,epsfig,verbatim}
\usepackage{amsmath,amsthm,amssymb,latexsym,graphics,textcomp,enumerate}
\usepackage[capitalize]{cleveref}
\usepackage{graphicx, psfrag}
\usepackage{imakeidx}
\usepackage{eucal} 
\usepackage[utf8]{inputenc}
\usepackage{xcolor}
\title[Marinus--Ptolemy and Delisle--Euler conical maps]{On the Marinus--Ptolemy and Delisle--Euler conical maps}
\author{Hideki Miyachi, Ken'ichi Ohshika and Athanase Papadopoulos \medskip}
\address{Hideki Miyachi,  
School of Mathematics and Physics,
College of Science and Engineering,
Kanazawa University,
Kakuma-machi, Kanazawa,
Ishikawa, 920-1192, Japan}
\email{miyachi@se.kanazawa-u.ac.jp}  
\address{Ken'ichi Ohshika,
Department of Mathematics,
Gakushuin University,
Mejiro, Toshima-ku, Tokyo, Japan}
  \email{ohshika@math.gakushuin.ac.jp}
\address{Athanase Papadopoulos,
Institut de Recherche Mathématique Avancée
(Université de Strasbourg et CNRS),
7 rue René Descartes,
67084 Strasbourg Cedex France and Centre for Interdisciplinary Mathematical Sciences, Institute of Science, Banaras Hindu University, Varanasi-221005, India. }
  \email{papadop@math.unistra.fr}

\numberwithin{equation}{section}

\newtheorem{theorem}{Theorem}[section]
\newtheorem{proposition}[theorem]{Proposition}

\theoremstyle{definition}

\theoremstyle{remark}

\begin{document}

\maketitle

\begin{abstract}

We examine connections between the mathematics behind  methods of drawing geographical maps due, on the one hand to Marinos and Ptolemy (1st-2nd c. CE) and on the other hand to Delisle and Euler (18th century). A recent work by the first two authors of this article shows that methods of Delisle and Euler for drawing geographical maps, which are improvements of methods of Marinus and Ptolemy, are best among a collection of geographical maps we term ``conical". This is an instance where after practitioners and craftsmen (here, geographers) have used a certain tool during several centuries, mathematicians prove that this tool is indeed optimal. Many connections among geography, astronomy and geometry are highlighted.
The fact that the Marinos--Ptolemy and the Delisle--Euler methods of drawing geographical
maps share many non-trivial properties is an important instance of historical continuity in mathematics.

\bigskip

\noindent {\bf Note} The final version of this paper will appear in \emph{Ga\d{n}ita Bh\=arat\=\i } (Indian Mathematics), the Bulletin of the Indian Society for History of Mathematics,
  
 \bigskip


\bigskip


 \end{abstract}
\section{Introduction}

In this paper, we highlight some relations between the works on cartography of, on the one hand, the two mathematicians, astronomers and geographers Marinos of Tyre (1st-2nd c. CE) and Claudius Ptolemy (2nd c. CE),  and on the other hand, the French geographer and astronomer Joseph-Nicolas Delisle (1688-1768) and the mathematician Leonhard Euler (1707-1783). It may be useful to recall that the latter was also a geographer and astronomer and that he left his mark in a number of other fields including mechanics, acoustics, optics and more. These four scientists, besides contributing to our knowledge in geography in the broad sense,  were involved in the practical construction of geographical maps. In fact, there are maps and methods of constructing maps known, quite rightly, as Marinos--Ptolemy and Delisle--Euler maps or methods. These maps also go by other names; indeed,  in cartography, unlike in mathematics, there is no uniform or  universally adopted terminology for objects. We shall see in this paper that Marinus--Ptolemy and Delisle--Euler  maps are constructed according to some common principles and meet the same needs.  In fact, our aim is to draw a parallel between the works of the two pairs of geographers, Marinos--Ptolemy and Delisle--Euler, and also to highlight the evolution of the methods of mapmaking known, perhaps heuristically, by Marinos and then Ptolemy in the 1st-2nd centuries CE, then studied more systematically by Delisle in the 18th century, until their theorization by Euler --- who was a collaborator of Delisle ---, using the newly developed methods of differential calculus and differential equations.  Finally, we shall highlight the fact that the Marinos--Ptolemy and Delisle--Euler methods, which have been in use for centuries and even millennia, have been shown, using modern computational tools, to have optimal properties among a family of maps that we call ``conical projections''. 

It is interesting to see that, once again, mathematical considerations confirm that a technique which has been used by practitioners (in the present case, Marinos, Ptolemy and Delisle, but there are also others) and which was studied in a purely mathematical framework by Euler, does indeed possess properties that make it an optimal method in a very precise sense which we shall make explicit. The fact that the Marinos--Ptolemy and the Delisle--Euler methods of drawing geographical maps share many non-trivial properties is an important instance of historical continuity in mathematics.

\noindent{\it Acknowledgement.}
We would like to thank the referee of this article for several interesting remarks.
 
\section{The problems of cartography}
  Let us start by briefly recalling the main problem in geography. At the same time, we shall point out an instance where some authors of
articles on history of mathematics sometimes repeat statements made by others without
checking the original sources, and errors spread out.

  In geography, one would like to represent on a Euclidean plane a country, usually considered to be a subset of the sphere, the latter representing the surface of the Earth.
  It is known since ancient times that such a representation cannot be faithful regarding distances. Indeed, there are theorems of spherical geometry, known since Greek Antiquity, that tell us that there is no mapping from a subset $\Omega$ of the sphere to the plane which preserves distances up to scale, provided $\Omega$ has non-empty interior. Let us state this as a proposition:

\begin{proposition}\label{prop:geography}
There is no map from any nonempty open subset of the sphere  to the Euclidean plane that  preserves distances up to scale.
\end{proposition}

This proposition follows from classical theorems in differential geometry. For instance, Gauss, in his famous memoir  \emph{Disquisitiones generales circa superficies curvas} (General investigations on curved surfaces) \cite{Gauss-English}, proved that the curvature of a surface is an intrinsic invariant. This implies that no open subset of the sphere can be sent isometrically to a subset of the plane, since the sphere has constant positive curvature and the plane has zero curvature.
But  the proof of \cref{prop:geography} does not need the full strength of Gauss's theorem.  Indeed, the result  follows easily from more classical results,  namely, propositions on spherical geometry that are contained in Menelaus's \emph{Spherics}, written in the 1st-2nd c. CE. Let us recall two of these propositions, each of which implies the result. 

The first proposition is Proposition 12 of the \emph{Spherics} which says the following:

 \begin{proposition}\label{prop:12}
 In any spherical triangle, the sum of the three angles is greater than
two right angles.
  \end{proposition}
  
  This implies \cref{prop:geography},  since it shows that no spherical triangle can be sent isometrically (even up to scale) into the Euclidean plane, and since any open subset of the sphere contains a spherical triangle.   
  
Another result from Menelaus's \emph{Spherics} that immediately implies \cref{prop:geography}  is Proposition 27, which says the following:

\begin{proposition}\label{prop:27}
  In any spherical triangle $ABC$,  let $D$ and $E$ be the midpoints of $AB$ and $BC$, and let $DE$ be the geodesic arc joining them. Then $DE>AC/2$.
\end{proposition}
   \begin{figure}[htbp]
\centering
 \includegraphics[width=0.7\linewidth]{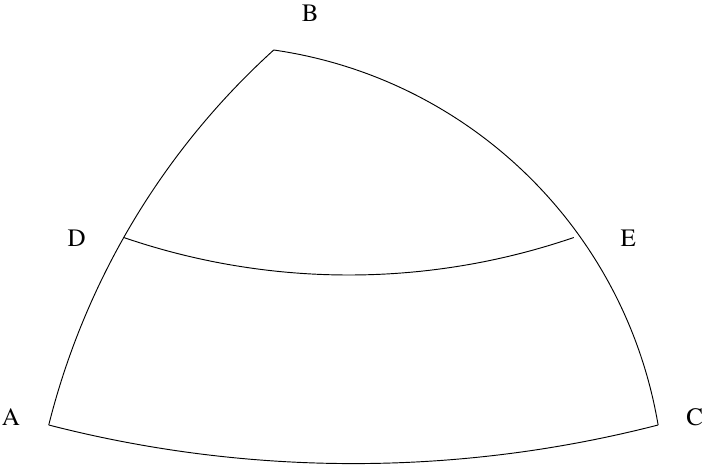}    \caption{\small Figure for \cref{prop:27}.} 
  \label{fig:M1} 
\end{figure}

 An English translation of Menelaus's \emph{Spherics} appeared recently \cite{RP}. The above two propositions are contained in p. 526 and 576 respectively. See also \cite{Menelaus-Monthly}.

Regarding this matter, we cannot help but note here that in the recent literature on mathematical cartography, Euler has been often misquoted. Indeed, he is credited with the fact of being the first to have shown that there is no map of the sphere to the plane that preserves distances up to scale. This erroneous claim is made in several relatively recent books and articles. We mention \cite{Feeman, Gray, Osserman, Rodriganez, Robinson} but there are many others.   The problem stems from a poor translation of Euler's memoir  \emph{De repraesentatione superficiei sphaericae super plano} \cite{Euler-De-represenatione} in which Euler proves that there is no ``perfect" map from a subset of the sphere to the Euclidean plane, and the error originates in a wrong interpretation of what Euler means by a ``perfect" map. 
   For instance,
R. Osserman, in his paper \emph{Mathematical mapping from Mercator to the millennium} \cite[p. 234]{Osserman}, referring to Euler’s memoir \cite{Euler-De-represenatione}, writes that the latter proved the following: ``It is impossible to make an exact scale map of any part of
a spherical surface", where an ``exact scale map" is defined by Osserman to be a map that preserves distances up to scale.  The problem is transmitted from one author to another  due to the fact that authors generally repeat information (or misinformation) they find in other articles, without checking the source. This problem and its effects is discussed in the introduction of the book \cite{CP}.  The fact that the translation of Euler's article on which several authors relied is not faithful was first pointed out by Charitos and Papadoperakis in their paper \cite{Ch-P}. The meaning of the adjective ``perfect" used by Euler is also reviewed in detail in the article \cite{Charitos}.    The translation problem was pointed out and discussed in the article \cite{Geography-GB} and in the opening pages of the book \cite{CP}. The faulty translation still circulates on the internet and problems of the kind we mentioned are still ongoing. Let us point out as an example the recent article that appeared in the \emph{Advances in Cartography and GIScience of the International Cartographic Association} \cite{Lapaine}, titled \emph{Map projection article on Wikipedia} \cite{Lapaine-Wikipedia}, in which the author makes a critique of the Wikipedia page on \emph{Map Projection}, but falls in the same trap. Indeed, he begins by quoting the following passage from the Wikipedia page in question:   
``[\ldots] However, Carl Friedrich Gauss's \emph{Theorema Egregium}
proved that a sphere's surface cannot be represented on a
plane without distortion," and he writes after this: ``The last sentence is correct, but
injustice has been done to a prominent mathematician,
physicist and cartographer Leonhard Euler (1707--1783).
He has authored more than 700 works etc.". Then, after quoting Euler's three articles on geography \cite{Euler-De-represenatione, Euler2, Euler-Delisle}, the author of \cite{Lapaine-Wikipedia}  writes: ``In the first of them, \emph{On mapping the sphere in the plane},
[Euler] gave the first formal proof of the impossibility of
mapping the sphere into the plane without any distortion
(Euler, 1777) [\ldots] C. F. Gauss was
born in 1777, and that same year, L. Euler published
proof of the impossibility of mapping the sphere in the
plane without any distortion. Therefore, the
championship belongs to Euler and should not be
overwhelmed." In fact, this article adds more to the existing confusion.

A lesson on writing on the history of mathematics emerges from all this.

In Sections \ref{s:Marinos} and \ref{s:Delisle}, we shall say a few words on our four geographers, Marinos, Ptolemy,  Delisle and Euler. We do not need to remind the reader who Ptolemy and Euler are, as their names are familiar to mathematicians, but we shall say a few words about them as geographers, an important attribute they had and which is not well known by these same mathematicians.
 
\section{Marinos and Ptolemy: Some biographical markers}\label{s:Marinos}

Marinos of Tyre (ca. 70-130 CE) was a mathematician and geographer, sometimes referred to as the father of mathematical geography and cartography. He was 30 years older than Ptolemy (ca. 100-168). 
Marinos was Phoenician and Ptolemy was Greek, and they both lived under Roman rule: The first  witnessed the reigns of Vespasian, Titus, Domitian, Trajan and Hadrian, whereas the second experienced the reigns of Trajan, Hadrian, Antoninus the Pious, Lucius Aurelius Verus and Marcus Aurelius. Ptolemy carried a Roman first name, Claudius. To Marinos is attributed a treatise, the \emph{Geography}, which did not survive but which is known to us because Ptolemy quotes it in his own \emph{Geography}.  As a matter of fact, all the information we have on Marinos comes from Ptolemy. The latter, in the  first 15 (out of a total of 24) chapters of Book I of his \emph{Geography}, heavily quotes Marinos' work, even if he often does so in order to criticize it. Book I of Ptolemy's \emph{Geography} is the theoretical part of his treatise. It is there that this author expounds  his methods of drawing geographical maps. We learn there that Marinos has developed a method of drawing geographical maps, and Ptolemy gives some details on this method. We shall say more on this in the next section. We also learn from the same work of Ptolemy that Marinos accurately noted a large number  of coordinates of places on Earth, by their latitude and longitude,  and that he himself used this data (with corrections).\footnote{According to Tannery, the technical term ``latitude" was first introduced by Ptolemy, see \cite[p. ii]{Tannery}.} 
Seven centuries later, the Arab historian and geographer al-Mas`\=ud\=\i \ (c. 996-956) cited Marinos's \emph{Geography} as a treatise to which he had access through an Arabic translation, but no manuscript of such a translation survives.
 
Ptolemy spent his life in Egypt, working mainly in Alexandria, a city founded by the Greeks and which, until the arrival of the Romans in 47 BCE, was the celebrated Greek city for what concerns science. He is credited with the most influential (in terms of duration) astronomical work of all times, mostly known by the Arabic name \emph{Almagest}. The original Greek name was \emph{Math\=ematik\=e S\'untaxis} (Mathematical composition); it was later renamed \emph{H\^e Meg\`al\^e S\'untaxis} (the Great composition), which became \emph{H\^e meg\'\i st\^e} (the Greatest) and which led, in Arabic, to the title \emph{Al-magis\d{t}\=\i }, transcribed in Latin as \emph{Almagestum} and which gave the name \emph{Almagest}, which lasts to this day. This work was translated into arabic by the Arab mathematician al-Khw\={a}rizm\={\i } (born c. 780). It was translated into Latin from the Arabic in the \textsc{xii}th century. It constitutes the most comprehensive surviving ancient treatise on astronomy. We refer the reader to the editions \cite{Ptolemy-Almagest} and \cite{Berggren}.

Two of Ptolemy's other important treatises survive: the \emph{Harmonics}, a treatise on music that remained  for several centuries the authoritative reference work on the theory and mathematical principles of music, and the one that interests us here, the \emph{Geography}, which constitutes a synthesis of all the geographical knowledge of the Greco-Roman world; see \cite{Ptolemy-Harmonics}. We also mention that Ptolemy wrote another geographical work, known in Latin
as the \emph{Planisphaerium} (Planisphere), which reached us only through Arabic translations, see \cite{Ptolemy-Planisphere}. In this work, Ptolemy studies the properties of the stereographic projection that he used for drawing maps of the celestial sphere, a projection which was especially useful in the conception of the Astrolabe.\footnote{The Astrolabe is an instrument that combines a representation of the celestial vault with a calendar. This instrument was highly developed by the Medieval Arab mathematicians for its usefulness in finding the direction of Mecca and the time for prayer, and this acted as a motivation for several Arab mathematicians of the Middle Ages to develop and improve its theory.}
We shall not linger over these works by Ptolemy, but we do mention that the three subjects, astronomy, geography and music, were closely related. The relationship between astronomy and music, theorized since the time of the Pythagoreans, was brilliantly developed by Johannes Kepler in his major work, \emph{Harmonices Mundi} (published in 1619) \cite{Kepler}. Regarding the relationship between astronomy and geography --- the subject that interests us here ---, it is important to know that since the epoch of the first great Greek geographers (Hipparchus, Eratosthenes and others), the coordinates of cities and the distances between distant points on the Earth's globe were measured using astronomical data.  We shall say more on this in \cref{s:Marinos--Ptolemy}. Let us also quote in this respect Ptolemy, who, by the end of Book II of the Almagest, after establishing the  tables of coordinates of the noteworthy stars, writes \cite[t. 1 p. 148]{Ptolemy-Almagest}:
\begin{quote}
\small This table of angles should be followed by the locations of the most famous cities, according
to their longitudes and latitudes, computed according to the celestial phenomena that are observed in each of these cities. But we shall treat distinctly this interesting subject which belongs to geography, and we shall use, for this purpose, the memoirs and the relations of
the authors who wrote on this matter. We shall indicate the place of each city by how many degrees, counted on the meridian, it is far from the equator; and in degrees, counted on the equator, the distance to the East or to the west of each meridian, to the one passing by
Alexandria [\ldots].
\end{quote}
Let us also quote the following passage from  Book VII of the \emph{Geography} \cite[33, p. 71]{Ptolemy-geo}: 
\begin{quote}\small
After I showed, at the beginning of the Mathematical Composition, how one could describe
on a sphere the known part of the Earth, and represent it on a planar surface in the most
similar and the most commensurably conformal manner to what we would see on the solid
sphere, it is necessary to present, as preliminaries of these expositions of all the inhabited
Earth, a graphical construction which covers all the parts that may be seen in each of these
tables. Thus, we shall add here this general representation, in its exact proportions.
\end{quote}
Finally, let us quote Nicolas Halma, the editor and author of the French translation of both Ptolemy's \emph{Almagest} and \emph{Geography}, who writes, in the Preface of \cite{Ptolemy-geo}:
\begin{quote}\small

 Ptolemy, descending from the observation of the heavens to the description of the Earth, here fulfils the commitment he made in his \emph{Treatise on Astronomy}, to provide a separate \emph{Geography} based on celestial observations. ``After teaching,'' says this Astronomer-geographer, ``the Method of finding the angles and arcs in each parallel, I would still have to speak of the position in longitude and latitude of the most remarkable cities of each country, to serve in the calculation of celestial phenomena seen from these cities. But we shall deal separately with this interesting subject, which belongs to Geography; and we shall rely, for this purpose, on the Memoirs and Relations of authors who have written on this subject. We shall mark by how many degrees, counted on its Meridian, the place of each city is distant from the equator; and in degrees counted on the equator, the Eastern or Western distance of each meridian, to that which passes through Alexandria; for it is to the meridian of this city that we compare those of the other points on the surface of the Earth."  \end{quote}

  There is a very large literature on the history of geography, written by geographers and by historians. The reader interested in a concise exposition of the history of Greek and Arabic geography may refer to the article \cite{Papa-Greek}, the first chapter of the book \cite{CP} which essentially concerns Eighteenth Century geography but where Ptolemy's work is put in perspective.

 \section{The Marinos--Ptolemy map}\label{s:Marinos--Ptolemy} 
 
 In this section, we shall essentially quote sections from Ptolemy's \emph{Geography} in which he describes Marinos's and his own methods of map drawing. Our goal is to highlight the similarities with the Euler and Delisle methods which we shall describe in \cref{s:Euler-proj}. We start by quoting Ptolemy from  \cite[Chap. XX, p. 82]{Berggren}:
 \begin{quote}\small

 Marinos paid considerable attention to this problem [of map drawing], and found fault with
absolutely all the [existing] methods of making plane maps. Nonetheless, he
himself turns out to have used the one that made the distances least proportionate.
He made the lines that represent the parallel and meridian circles all
straight lines, and also made the lines for the meridians parallel to one another,
just as most [mapmakers] have done; but he kept only the parallel through
Rhodes proportionate to the meridian in accordance with the approximate ratio
of 5:4 that applies to corresponding arcs on the sphere (that is, the ratio of the
great circle to the parallel that is 36$^{\mathrm{o}}$ from the equator), giving no further thought
to the other [parallels], neither for proper proportionality nor for a spherical
appearance. Now when the line of sight is initially directed at the middle of the northern quadrant of the sphere, in which most of the oikoumene is mapped, the meridians can give an illusion of straight lines when, by revolving [the globe or the eye] from side to side, each [meridian] stands directly opposite [the eye] and its plane falls through the apex of the sight.   The parallels do not do so, however, because of the oblique position of the north pole [with respect to the viewer]; rather, they clearly give an appearance of circular segments bulging to the south.
 \end{quote}
 
We also quote the following passage, from the next chapter  \cite[Chap. XXI, p. 83]{Berggren}, where Ptolemy describes his own method:
  \begin{quote}\small
 [\ldots] For these reasons it would be well to keep the lines representing the meridians
straight, but [to have] those that represent the parallels as circular segments
described about one and the same center, from which (imagined as the north
pole) one will have to draw the meridian lines. In this way, above all, a semblance
of the spherical surface will be retained in both its actual disposition
and its visual effect, with the meridians still remaining untilted with respect to
the parallels and still intersecting at that common pole. Since it is impossible to
preserve for all the parallels their proportionality on the sphere, it would be
adequate (1) to keep this [proportionality] for the parallel through Thule and
the equator (so that the sides that enclose our [oikoumene's] latitudinal dimension
will be in proper proportion to their true magnitudes), and (2) to divide the
parallel that is to be drawn through Rhodes (on which most of the investigations
of the longitudinal distances have been made) in proportion to the meridian,
as Marinos does, that is, in the approximate ratio of similar arcs of 5:4 (so
that the more familiar longitudinal dimension of the oikoumene is in proper
proportion to the latitudinal dimension).  
 \end{quote}

Thus, Ptolemy demands that the following are satisfied:

\begin{enumerate}
\item  \label{enum:1}The images of the meridians are straight lines that intersect at a common point.

\item   \label{enum:2}The images of the parallels are arcs of circles that have a common centre. This point is the image of the North pole, and it is also the intersection point of the images of the meridians.

\item  \label{enum:3} Distances are preserved up to a common scale for the parallels passing through the island of Thule and the equator.

\item  \label{enum:4} Distances are preserved up to scale for the image of the parallel passing through the island of Rhodes,  and this scale is $\frac{4}{5}$ if we consider that for the two  parallels mentioned in (3) above, the scale is 1.

\end{enumerate}

The reason of the choice of the two specific parallels in (3) is that they roughly correspond  to the upper and lower borders of the Known Inhabited World (the ``Oikoumene"). The reason of the choice of the parallel through Rhodes is that this island was, in the 1st-2nd centuries, an important Astronomical centre, similar to Alexandria. Incidentally, it was the place where the astronomer, geographer and mathematician Hipparchus (c. 190-c. 120 BCE) used to conduct his observations.\footnote{Pierre Simon de Laplace, in a passage of his major work \emph{Exposition du système du monde} (Exposition of the system of the world), explaining the importance of geography in the sciences, attributes to Hipparchus a method of computing coordinates by latitude and longitude using astronomical data, in particular, the eclipses of the moon \cite[vol. 2, p. 505]{Laplace}. The French 19th century mathematician, astronomer and historian of astronomy Jean-Baptiste Joseph Delambre called  Hipparchus ``the real father of astronomy" \cite[vol. 1, 2.6]{Delambre}.}

\begin{figure}[htbp]
\centering
 \includegraphics[width=\linewidth]{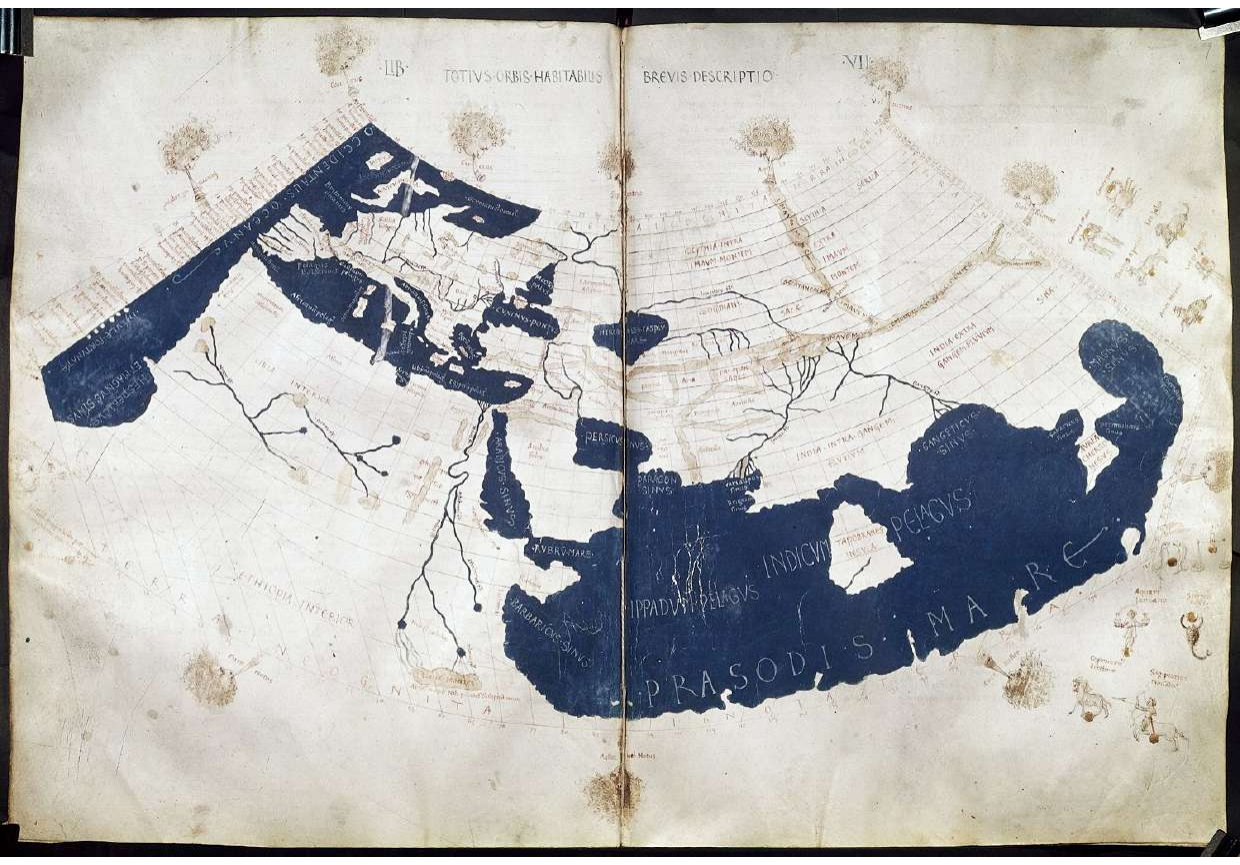}    \caption{\small A 15th century reconstruction of Ptolemy's world map. The British Library, Harley MS 7182.} 
  \label{fig:PtolemyWorldMap} 
\end{figure}

In \cref{fig:PtolemyWorldMap} we have reproduced a 15th century reconstruction of Ptolemy's world map.\footnote{This map is based on Jacobus Angelus's 1406 Latin translation of Maximus Planudes's late-13th century rediscovered Greek manuscripts of Ptolemy. }

Now we come to Delisle and Euler.

\section{Delisle, Euler and geography}\label{s:Delisle}
  
   Joseph-Nicolas Delisle (1688-1768) was a French physicist and astronomer,  member of the Paris Royal Academy of Sciences. He belonged to a family of famous  astronomers and geographers, and he was among the first scientists invited to Saint Petersburg by the Tsar Peter the Great at the foundation of the Russian Academy of Sciences. Delisle was appointed the director of the astronomical department at this Academy. During his long stay in Saint Petersburg (1727 to 1748), he founded the first Russian astronomical school, a school that became, a few decades later, one of the most renowned in the world. He was also in charge of the geography department at the Academy, where his main task was to conduct distance measurements and to collect geographical information that can be used for the drawing of reliable maps of the mostly uncharted huge Russian Empire. The two departments to which Delisle belonged, geography and astronomy, were part of the ``class" of mathematics,\footnote{The Russian Academy had eleven professorships, distributed into three ``classes" (departments): mathematics (in the 18th-century sense, which also included astronomy, geography, mechanics and naval science), physical sciences (which included  chemistry, botany and anatomy) and humanities
(which included rhetoric, history and law).} and Delisle, in his work, needed the help of mathematicians. He first collaborated for several years with the mathematician Jakob Hermann who, like him, was a founding member of the Saint Petersburg Academy, but in 1731, Herman returned to his home city, Basel (which, incidentally,  is also the home city of Euler). Delisle needed a new assistant, and the choice fell on Euler. 

Euler was the ideal collaborator for Delisle. He had already worked on astronomy, notably on the motion of the moon, and later, on the motion of comets. His first book on astronomy is titled 
\emph{Theoria motuum planetarum et cometarum} (Theory of the motions of planets and comets) \cite{E66}. 
It is likely that Euler had been influenced by Delisle in his choice of the astronomical problems he had worked on. He wrote several memoirs on astronomy, in particular on the determination of  orbits of comets and planets, on methods of calculating the parallax of the sun, and on celestial mechanics, a topic for which he won several prizes from the French Academy of Sciences, on questions on which he was competing with several other French and British scientists. Among his  works on astronomy we mention the book  \emph{Theoria motus lunae exhibens omnes eius inaequalitates} (Theory of the motion of the moon which exhibits all its irregularities) \cite{E187}, written  in 1751, in which he develops the so-called first Euler lunar theory, which includes a method of approximation for the three-body problem. We also mention the  \emph{Theoria motuum lunae} (Theory of lunar motion),  published in 1772 and  which contains the so-called Euler second theory of lunar motion.\footnote{One should note that the determination of the lunar motion was one of the most challenging astronomical problems since Ancient times, and it remained so until  the epoch of computers. As a matter of fact, this problem was the main object of the so-called three-body problem. Newton, before Euler, tried to deal with this problem in his \emph{Principia}, and the same problem was the  main object of study in Poincaré's famous work, \emph{Les méthodes nouvelles de la
mécanique céleste}, written in three parts and published in 1892, 1893 and 1899 \cite{Poincare}.}

Euler, during several years, assisted Delisle in his work on geography, a task in which he was thoroughly involved. 
In 1738, he lost the use of his right eye  and he attributed this fact to the excessive stress due to his work on the examination of geographical maps.\footnote{Euler wrote to Christian Goldbach, on August 21st, 1740 \cite[Vol. 1, p. 163, English translation p. 670-671]{Euler-Goldbach}: ``Geography is fatal to me. As you know, Sir, I have lost an eye working on it; and just now I nearly risked the same thing again. This morning I was sent a lot of maps to examine, and at once I felt the repeated attacks. For as this work constrains one to survey a large area at the same time, it affects the eyesight much more violently than simple reading or writing. I therefore most humbly request you, Sir, to be so good as to persuade the President by a forceful intervention that I should be graciously exempted from this chore, which not only keeps me from my ordinary tasks, but also may easily disable me once and for all."} 

Delisle defined a triangulation of Russia, for his project of constructing maps of the various regions of the Russian Empire. Let us recall that the so-called triangulation method allows the measurement of long geodesic arcs, or of distances between two points on the surface of the Earth  that are far from each other (at the scale of distances between two cities). In fact, this method reduces the measurement of such an arc or distance to a measurement of angles. This is due to the spherical shape of the Earth. Mathematicians can easily understand the theory behind this:  Suppose we want to measure the length of a long geodesic arc on the surface of the Earth, supposed to be spherical. The geographer starts by defining a  triangulation consisting of a collection of geodesic triangles that are located on both sides of this geodesic arc.\footnote{A geodesic triangle on a surface is the union of three
geodesic arcs that meet pairwise at a point called a vertex of the triangle.} The important fact used is that a spherical triangle  (unlike a Euclidean triangle) is completely determined by its three angles. As we already mentioned, astronomical observations were used to measure the angles at the vertices of these triangles. Then, from the spherical trigonometric formulae, one can find the edge lengths. Once the various lengths of the edges of these triangles are determined, one can find the length of the initial arc. Triangulations were used in geography since Antiquity. This is reported on by Euclid,
Heron and others, see e.g. \cite[p. 845 and 893]{Neugebauer}. See also the comments in \cite[p. ii]{Ptolemy-geo}. It is interesting to know that Gauss\footnote{In particular, Gauss, in 1818, was getting funded by the King of Hanover  for a project of triangulating the Hanover region, of which Göttingen (the city were Gauss used to live) was part. The project took 30 years for its realization.} worked on the development of the theory of triangulations of surfaces in geometry, after using extensively the method of triangulations in geodesy. Delisle also used triangulations to conduct measurements that supported his view on the shape of the Earth.\footnote{Delisle shared the Newtonian conviction  that the Earth is flattened at the poles (a hypothesis which was later found to be the right one), rather than the Cassinian one, saying that the Earth is flattened at the equator. Delisle, the Cassinis and  other geographers used triangulations to calculate lengths of arcs of meridians, to confirm their theories on the flattening of the Earth. We also note that the discovery of the fact that the Earth is spheroidal and not spherical had impacts on works on geography. For instance, Euler, wrote a memoir titled \emph{\'Elémens de la trigonométrie sphéroïdique tirés de la méthode des plus grands et plus petits} (Elements of spheroidal trigonometry extracted from the method of maxima and minima) \cite{Euler-Spheroid}, in which the applications of the theory he develops concern distances between places on the Earth, assumed to be spheroidal.  The interested reader may refer to the survey and the discussion in \cite{Papa-Maupertuis} and \cite{Clairaut}.}

\begin{figure}[htbp]
\includegraphics[width=12.5cm]{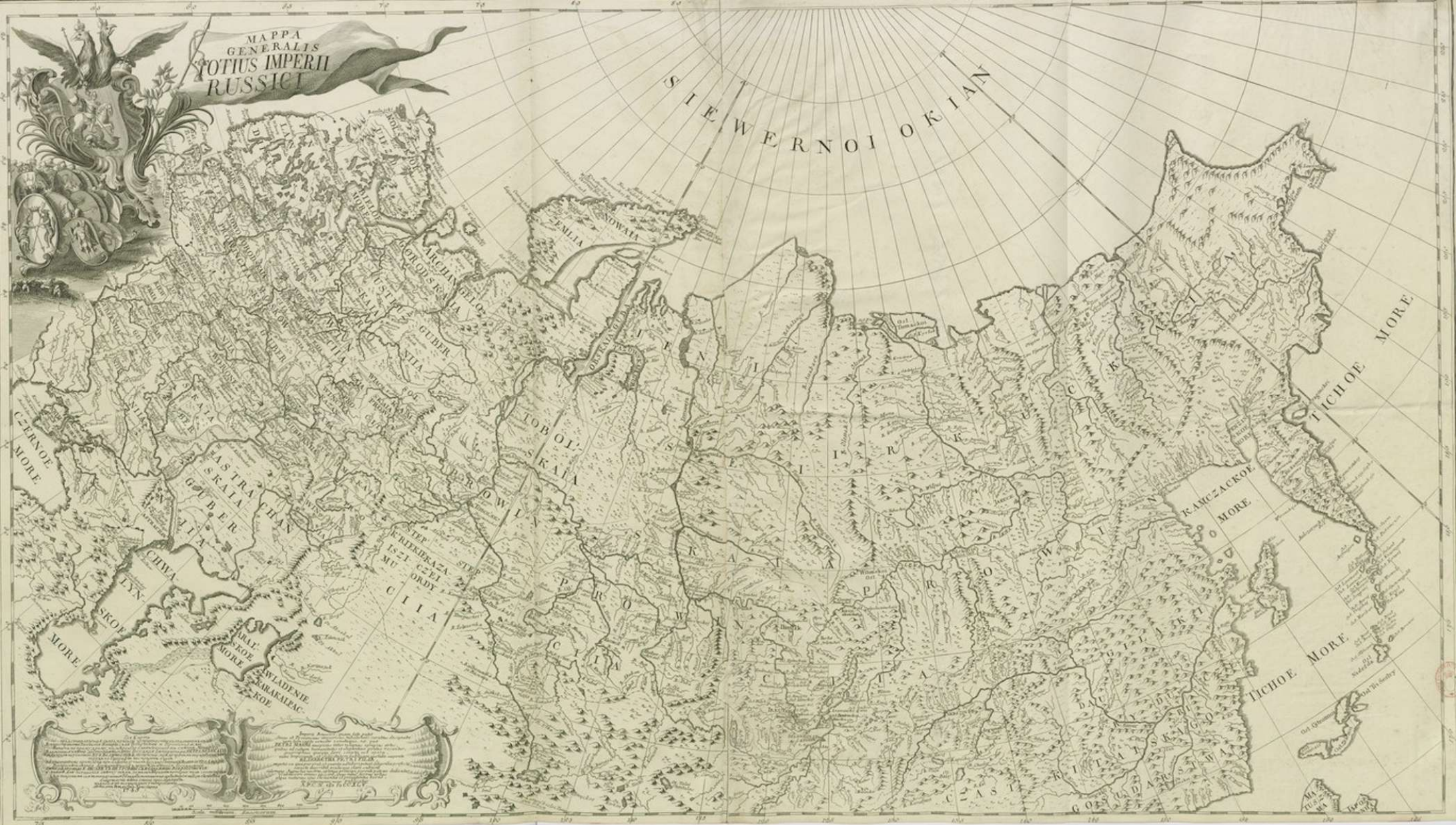} 
\caption{\small A map of the Russian Empire: this is the last map of the \emph{Atlas Russicus} (Saint Petersburg, 1745), and it is drawn using Delisle's method. Bibliothèque
Nationale de France, Département Cartes et Plans.}
 \label{fig:Empire} 
\end{figure}

In 1740, after some serious disagreements between Delisle and the administration of the Academy, Euler became the official leader of the project of establishing maps of the Russian Empire. He worked on this project in association with Gottfried Heinsius (1709-1769), a German mathematician and astronomer who, like Euler and Delisle, had been hired at the Russian Academy. At the same time, Euler  practically became in charge of the publication of the so-called \emph{Atlas Russicus} (Russian Atlas),\footnote{The complete title is 
\emph{Atlas russicus mappa una generali et undeviginti specialibus vastissimum imperium russicum
cum adiacentibus regionibus secundum leges geographicas et recentissimas obsevariones
delineatum exhibens. Cura et opera Academiae Imperialis Scientiarum Petropolitanae} (Russian atlas containing a general
map and nineteen particular maps of the whole Russian Empire and of the limiting countries,
constructed according to the general rules of geography and the most recent observations).} edited
by the Academy and which was eventually published under Delisle's name in 1745 \cite{Euler-Atlas-Russicus}. In \cref{fig:Empire} we have reproduced a map of the Russian Empire extracted from Euler's Russian Atlas, drawn according to Delisle's method. Eight years later, another atlas was published under the auspices of the Prussian Academy of Sciences in Berlin, where Euler was settled and under his editorship. This atlas is known under the name \emph{Atlas Geographicus}\footnote{The complete title is \emph{Atlas Geographicus omnes orbis terrarum regiones in XLIV
tabulis exhibens} (Geographical atlas representing in 44 tables all the regions of the
Earth).} \cite{Euler-Atlas-geographicus}. The maps in this atlas represent  large portions of the world. Euler wrote the preface and  commented on each of the maps it contains. 

\begin{figure}[htbp]
\centering
 \includegraphics[width=\linewidth]{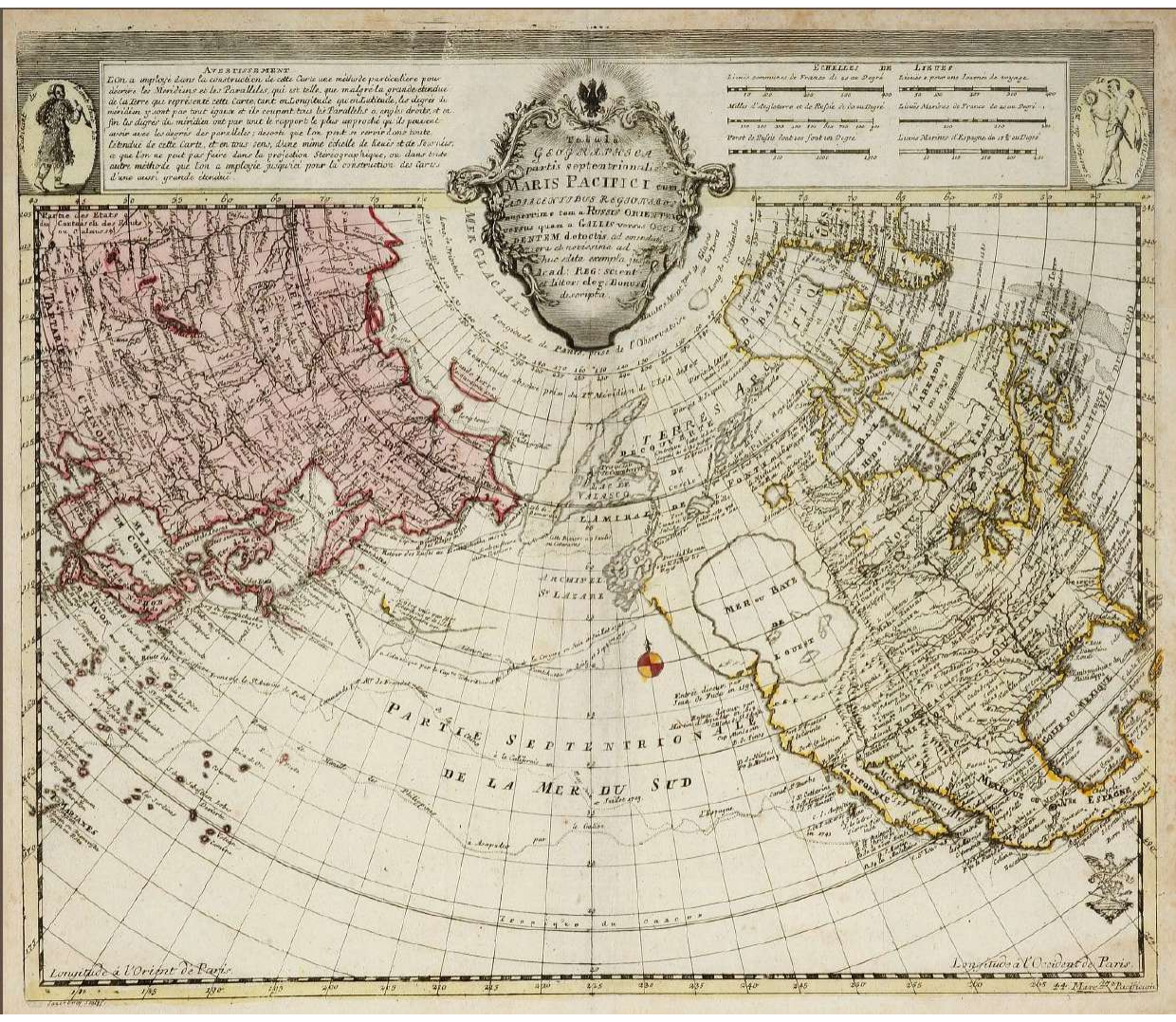}    \caption{\small Map of the Northern Pacific, representing the Eastern part of Asia and the Northern part of America, from Euler's \emph{Atlas Geographicus} (Berlin, 1753) and it is drawn using Delisle's method.} 
  \label{fig:Euler-Map41} 
\end{figure}

Delisle remains one of the scientists who made a substantial contribution to geography, through an extensive use of astronomical observations. His work on geography, together with Euler's work in this field, are introduced and commented on in the volume \cite{CP}.  The same volume also contains a translation of   the part of the preface of Euler’s
Berlin atlas in which the latter comments on  the map reproduced here in \cref{fig:Euler-Map41}, which is the last one in the
series and which is drawn using Delisle’s method.

 In 1777, Euler wrote three memoirs on geography,  \emph{De repraesentatione superficiei sphaericae super
plano} (On the representation of the spherical surface on the plane) \cite{Euler-De-represenatione}, \emph{De proiectione geographica superficiei
sphaericae} (On the geographical projection of the
surface of the sphere) \cite{Euler2} and \emph{De proiectione geographica De Lisliana in mappa
generali imperii russici usitata} (On Delisle’s geographical projection used for a general map of the Russian
Empire), see \cite{Euler-Delisle}. 
  In the next section, we shall say more about Euler's work on geography.

\section{The Euler-Delisle projection}\label{s:Euler-proj}

Euler, in the preface of the
 atlas published by the Prussian Academy \cite{Euler-Atlas-geographicus},  gave some indications on the methods used in the drawing of each of the 44 maps that constitute it, see \cite{Atlas-Intro}. Several kinds
of projections are used in drawing these maps, depending on the shape of the country represented. In all of them, the meridians are
represented by straight lines and the parallels are represented by either straight lines or circles that are perpendicular or almost perpendicular to the meridians.   We have reproduced in \cref{fig:Euler-Map41} the map of the Eastern
part of Asia and the Northern part of America, which is drawn using Delisle’s method. Euler's comments on this map are translated in \cite[p. 134]{CP}. He starts with the sentence: ``In the description of this map, we have kept the method which the famous Mr. Delisle used,
which seems to us the most appropriate for a good representation of these Northern regions.
In the general map of the Russian Empire published by the Saint Petersburg Academy, we
had used this same method, which at first glance may not be satisfying; we shall explain it
here in a few words." We summarize some properties of this projection, as Euler presents it.   
  \begin{enumerate}
         
\item The region drawn is contained between elevation 45$^{\mathrm{o}}$ to 68$^{\mathrm{o}}$ from the pole.\footnote{Euler calls ``elevation from the pole" the co-altitude, that is the difference from 90$^{\mathrm{o}}$  of the latitude. Note that the figures given here do not correspond to the two parallels that bound the Russian Empire (at the time of Euler). Euler, at the beginning of  \S 7 of  the paper quoted, takes 45$^{\mathrm{o}}$ to 68$^{\mathrm{o}}$ to be the geographical co-latitudes of a random  meridian  arc passing through the Russian empire (he writes: ``let AB [in the map] be an arc of an arbitrary meridian passing
through the Russian Empire, having in $A$ the southernmost extremity and in $B$ the
northernmost extremity "). A few lines after this, he writes, to justify this choice: ``But actually,
as in the general map of the Russian Empire, beyond the 70th degree of latitude
[only] unimportant regions are taken into account in order to be represented [\ldots]"}

\item On the map, the meridians  intersect at a common point, which is not the North Pole but is 7$^{\mathrm{o}}$ farther than this pole. 

\item The parallels are  circles that are equally distant.

 \item Among the meridians that are drawn, any two that are distant from each other by one degree converge in such a manner that under any two latitudes,  the ratio the degrees of longitude and of latitude are the same as in reality.
\end{enumerate}

Euler writes that in this way, the distances measured on the map  differ from the true ones by an amount  that can hardly be seen. He adds that if we wanted to draw using the same method regions that are close to the Pole or to the Equator, the error would be very large and that Delisle was aware of that. He notes that we should not regard as a defect  the fact that the centre at which all the meridians meet is not the Pole and that we should not be surprised by the fact that on this map the parallels that form the semi-circles do not occupy 180$^{\mathrm{o}}$ in longitude, but much more, sometimes even up to 250$^{\mathrm{o}}$: this shows, he says, that this map does not suffer in a large Latitude extent.

In \S 5 of the memoir  \emph{On Delisle's geographical projection used in the general map of the Russian empire} \cite{Euler-Delisle},
Euler states four properties that are ideally required from a geographical map:
\begin{enumerate}

\item The images of the meridians are straight lines; 
 
\item the degrees of latitudes do not change along meridians;
 
\item the images of the parallels meet the images of the meridians at right angles; 

\item at each point of the map, the ratio of the degree on the parallel to the degree on the meridian is the same as on the sphere.
\end{enumerate}

He notes that since these conditions cannot be achieved simultaneously, because of the last one, it is required, instead of this condition,  that the deviation of the degree of latitude to the degree of longitude at each point from the true ratio be as small as possible (ideally, this error should be unnoticeable). He  recalls then that Delisle was in charge of constructing such maps, and that in doing so he started by requiring that the ratio between latitude and longitude be exact on two noteworthy parallels, and that the question becomes that of choosing these two  parallels. 
Starting from \S 7 of his memoir,  Euler develops some mathematical tools that show how to construct such a map. This reminds us of Ptolemy's requirement that distances be faithful on the parallel of Rhodes.
Euler starts by choosing a meridian passing through the Russian Empire, considered  as the principal meridian, and then 
 he gives a method for drawing the other meridians, a method that  allows the construction of circles representing the parallels, such that for the  region comprised between two latitudes given in advance, the ratio of the degree of latitude to  the degree of longitude is faithfully represented. 
In \S 10--16, he gives the mathematical details showing how far this representation differs from reality at the extreme points he started with, and in \S 17--23, he makes the actual computations in the  special case where the map is that of the Russian Empire.
In the article \cite{Charitos-Papa}, Delisle's method, after Euler, is commented on, from the mathematical point of view.
 
  For a review of and a commentary on the general topic of map drawing in the works of Euler and Delisle, the reader can also refer to the article \cite{Geography-GB} and to the book on mathematical geography in the eighteenth century \cite{CP}, in particular the chapter titled \emph{Euler, Delisle and Cartography} \cite{Delisle-AP}.

 The next section contains a discussion on the name of Euler's projection and that of Ptolemy.

\section{About terminology in geography: conical and other projections} \label{s:termiology}

 The list \emph{Terminology of map projections}, compiled by D. H. Maling and presented at the Third International Conference on Cartography in Amsterdam (17-22 April 1967) \cite{Maling}, includes 234 entries, sometimes with sub-entries (234a, 234b, etc.), indicating the terminology for methods of map drawing and names of these maps in English, French, German and Russian. Several names indicate a conical projection; for example, we read (for names in English), at No. 48: \emph{Conical equidistant projection with one standard parallel} (Ptolemy), at No. 49: \emph{Conical equidistant projection with two standard parallels (Delisle)}, at No. 50: \emph{Conical equidistant projection with two standard parallels (Euler)}, at No. 54: \emph{Conical equidistant projection with two standard parallels (Murdoch I)}, and at No. 55:  \emph{Conical equidistant projection with two standard parallels (Votkovsky I)}.

There is a large number of  variations of  conic map projections in the various geographical articles on geography or history of geography, and in catalogues of geographical maps. 
Terminology in cartography is very far from being uniform and universally accepted. The Ptolemy--Marinus and the Delisle--Euler maps are sometimes called conical maps. In this section, we shall say in what sense we think that this name is deserved.

Let us start with the fact that the most basic definition of a conical geographic projection is a projection in the geometric sense, from a point or parallel to some direction, of the sphere onto a cone. The cone may touch the sphere in one circle, or it may intersect it in two circles, and its apex may or may not be on the sphere. Then, when the cone is developed on a Euclidean plane,\footnote{The cone is a developable surface, that is, it can be developed on a plane.} the image of the land we want to map becomes an image on the plane. This easily leads to a rigorous mathematical definition of a conical map, however restrictive. A mathematician will also be able to formalize the definition of a family of conical projections, which give as limits, on the one hand, a cylindrical projection and, on the other hand, a projection on a plane, the latter being what we usually call a stereographic projection.

The Ptolemy--Marinus and the Delisle--Euler projections are not conical in the above restricted sense, but both can be termed conical, for several reasons that we explain now.

Let us come back to what Ptolemy writes in his \emph{Geography}. An excerpt  we quoted in \cref{s:Marinos--Ptolemy} says: ``Now when the line of sight is initially directed at the middle of the northern quadrant of the sphere, in which most of the oikoumene is mapped, the meridians can give an illusion of straight lines when, by revolving [the globe or the eye] from side to side, each [meridian] stands directly opposite [the eye] and its plane falls through the apex of the sight.   The parallels do not do so, however, because of the oblique position of the north pole [with respect to the viewer]; rather, they clearly give an appearance of circular segments bulging to the south."

Clearly, Ptolemy is thinking of the following:  We look at the sphere, on which the  world map is drawn,\footnote{In fact, Ptolemy, before the passage quoted below, discussed a way of representing the Earth on a sphere, and it is to this latter representation that he refers.} by putting the eye in front of a meridian in such a way that we see this meridian as a straight line. Ptolemy then says that when the Earth is rotated around the axis passing through the pole and hidden by the meridian in front of the eye, the other meridians appear one after the other as straight lines, having a common intersection point. Seen in this way, there is no difference between the projection he describes and a conic projection in the above mathematical sense.

\emph{Let us note in passing that Ptolemy, without meaning it, is assuming that the Earth is rotating around an axis!}

\begin{figure}[htbp]
\centering
 \includegraphics[width=0.3\linewidth]{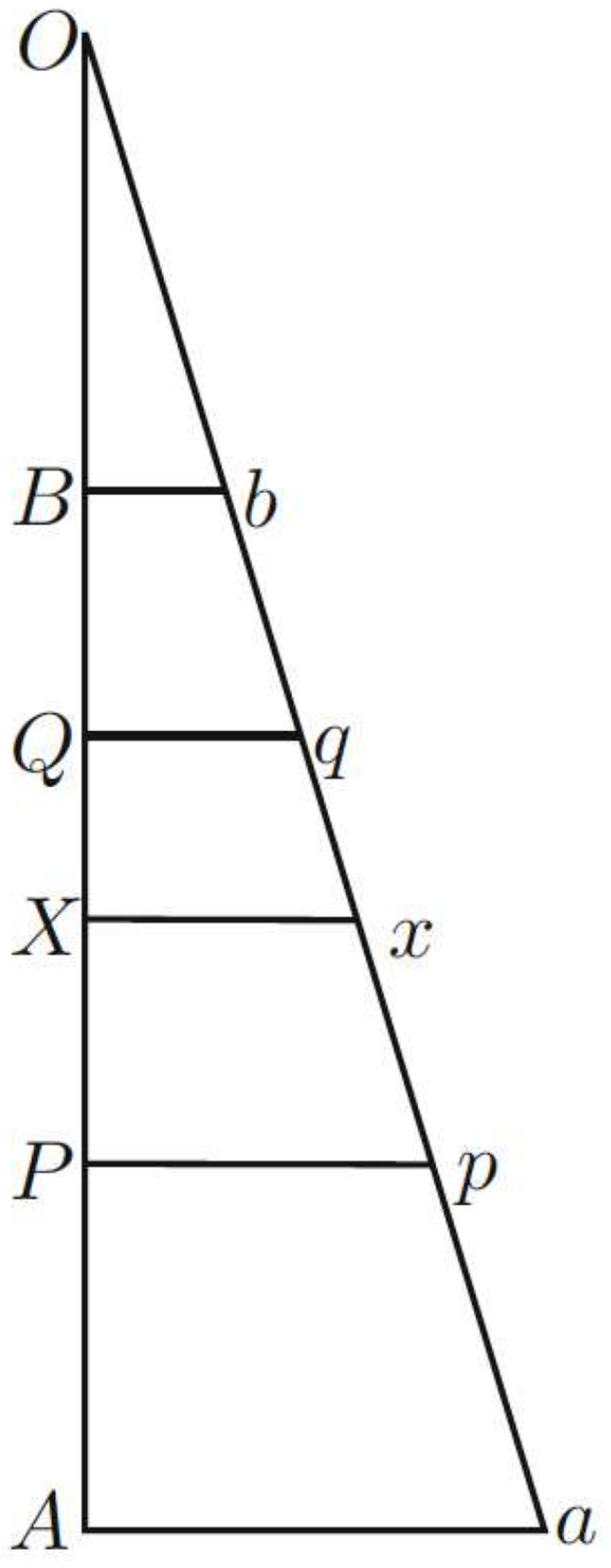}    \caption{\small A figure extracted from Euler's memoir \cite{Euler-Delisle} in which he explains Delisle's projection. The reader may refer to that memoir or to the paper \cite{Charitos-Papa} for the explanations.} 
  \label{f:Euler} 
\end{figure}

As for the Delisle--Euler projection, it also has all the characteristics of a conic projection in the above sense. Euler described the construction of the images of the meridians, one by one, and the construction can be done easily and naturally, like that of Marinos--Ptolemy, by passing through a projection on a cone. We have reproduced in \cref{f:Euler} a figure from Euler's memoir which accompanies his explanations for this construction.

Before concluding, we need to recall a few words on the theory of quasiconformal mappings and the Tissot indicatrix.
    
\section{The Tissot indicatrix of a geographical map} 

 The Tissot indicatrix is a device introduced by the French mathematician and geographer Nicolas-Auguste Tissot (1824--1897)  who studied extensively, from the mathematical viewpoint, the quasiconformal distortion --- that is, a measure of the lack of conformality, or of angle-preservation --- of mappings from the sphere on the Euclidean plane, with a view on applications to the drawing of geographical maps. He also developed a theory of quasiconformal distortion of mappings between general surfaces. His ideas are at the basis of the work on quasiconformal mappings that was done several decades after him by Gr\"otzsch, Lavrentieff, Ahlfors and Teichm\"uller.\footnote{Teichm\"uller writes in in §33  of his groundbreaking paper  \emph{Extremale quasikonforme Abbildungen und quadratische Differentiale} (Extremal quasiconformal mappings and quadratic differentials) \cite{T-1940}, that the origin of the  notion of quasiconformal mappings lies in the work of Tissot on the drawing of geographical maps.}  
 Let us be more precise.

\begin{figure}[htbp]
\centering
 \includegraphics[width=0.7\linewidth]{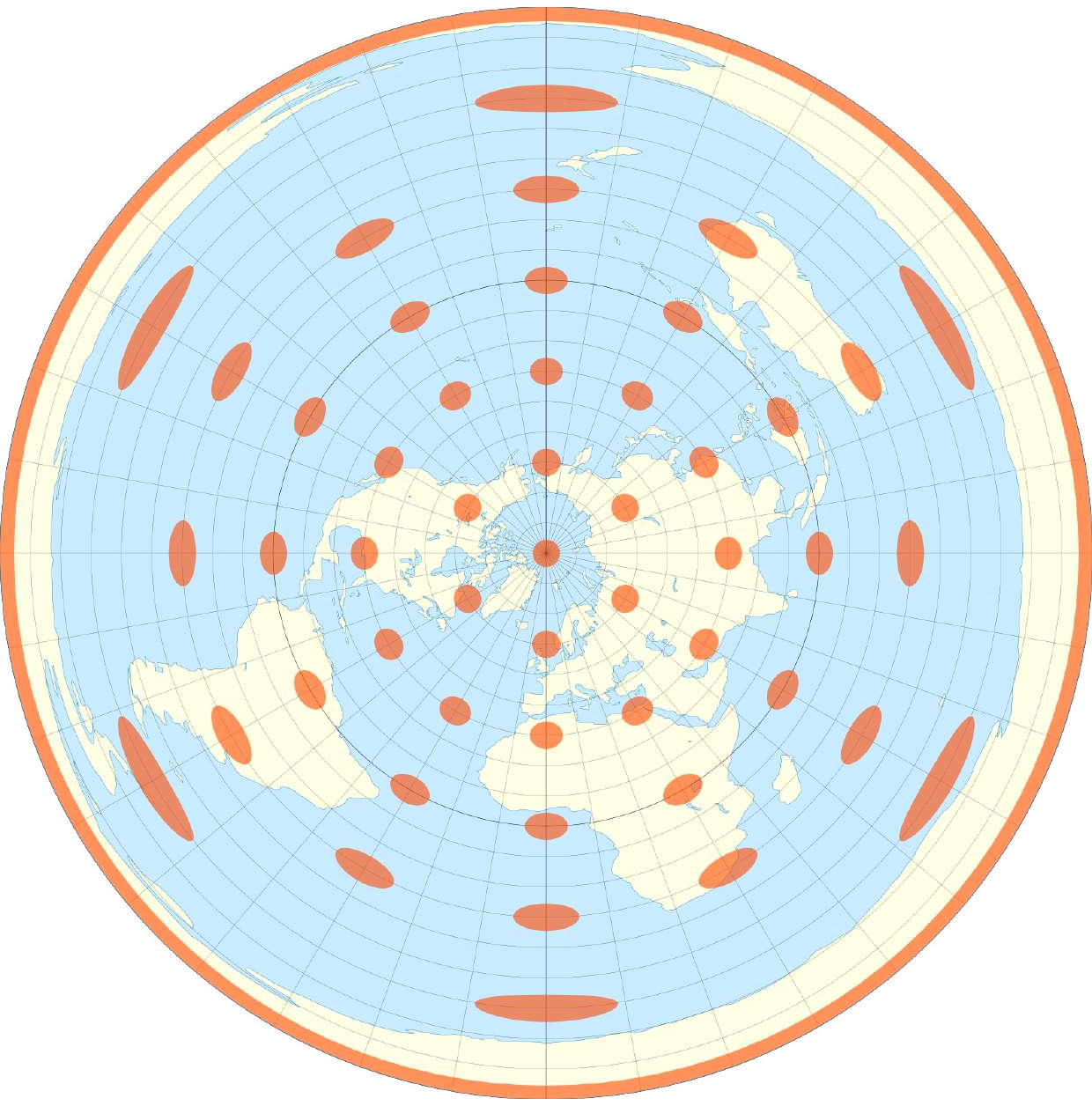}    \caption{\small The Tissot indicatrix for the azimuthal equidistant geographic projection.} 
  \label{f:azimuthal} 
\end{figure}

\begin{figure}[htbp]
\centering
 \includegraphics[width=0.8\linewidth]{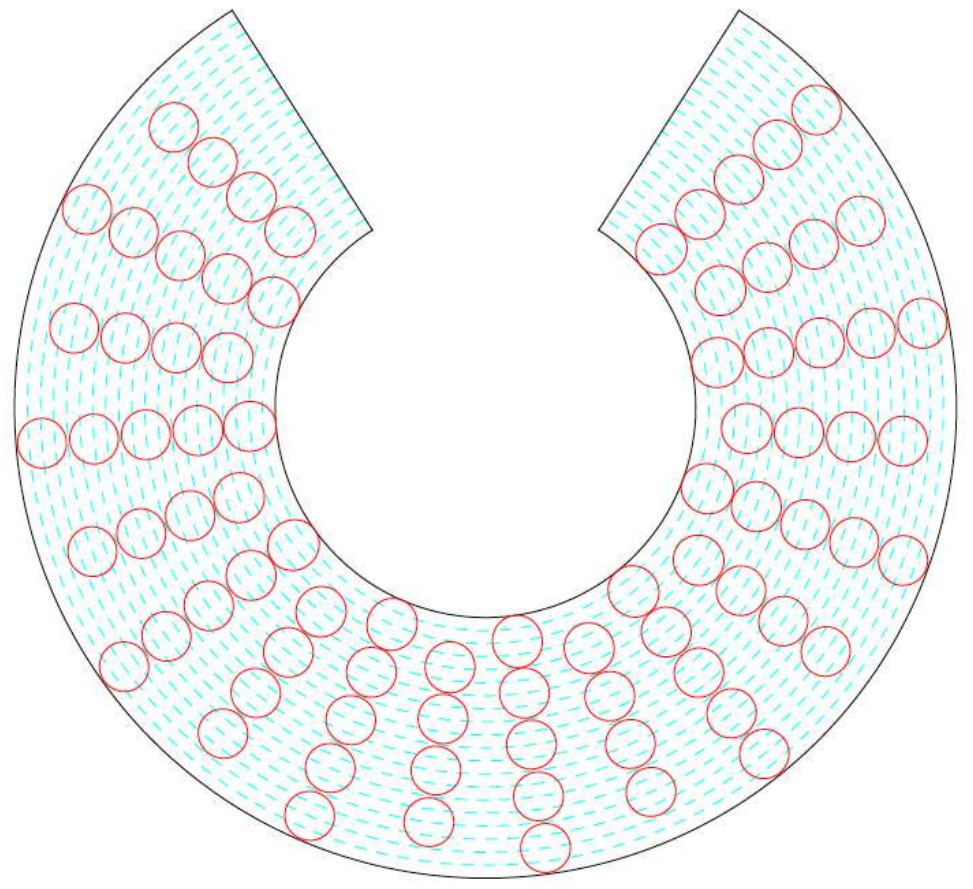}    \caption{\small The Tissot indicatrix for the Delisle--Euler map, reproduced with permission from the paper \cite{MO}. The reader will notice that the ellipses are almost circles of the same size, which shows that this map is almost angle-preserving.} 
  \label{fig:f3_round} 
\end{figure}

Tissot introduced a device that geographers call \emph{indicating ellipse} or \emph{Tissot indicatrix}. This device has been used by geographers until the middle of the twentieth century. It is a field of ellipses drawn on the map,  each ellipse representing the image by the projection (assumed to be differentiable) of an infinitesimal circle\footnote{Assuming, to simplify the matter, that   the map is differentiable, the expression ``infinitesimal circle" means  a circle centred at the origin in the tangent space at a point. In practice, it is a circle with a ``tiny radius" on the surface. In the art of geographical map drawing, these circles, on the domain surfaces, are all supposed to have equal sizes, so that the collection of relative sizes of the image ellipses becomes also a meaningful quantity.} on the sphere. Examples of Tissot indicatrices for geographical maps are given in \cref{f:azimuthal}  and   \cref{fig:f3_round}. The first figure represents the azimuthal equidistant geographic projection.  In  \cref{fig:f3_round}, we do not see the geographical map, but this figure represents the field of Tissot ellipses for the Delisle--Euler map, as we shall see below.

From a differential-geometric point of view, the Tissot indicatrix is an information on the metric tensor obtained by pushing forward the metric of the sphere by the geographical projection.   In the modern development of quasiconformal theory by the four authors we mentioned above, a significant parameter of a map is its quasiconformal dilatation at a point. This is the ratio of the major axis to the minor axis of the infinitesimal ellipse which is the image of an infinitesimal circle by the map.
The Tissot indicatrix provides more information than the quasiconformal dilatation, since it keeps track of (1) the direction of the great and small axes of the infinitesimal ellipse, and (2) the size of this ellipse, compared to that of the infinitesimal circle of which it is the image. 
 See \cite{Tissot1881} for Tissot's work. For a review of this work, the reader may refer to the paper \cite{Papa-Tissot}.

\section{Comparison between conical projections: dilatation and quasiconformal distortion} \label{s:qc}

In the recent paper \cite{MO} by the first two authors of the present paper, the Delisle--Euler map has been compared with several other maps in the category of conical maps. Using modern computing tools, these authors obtained numerical values for global metrical and quasiconformal distortions of these mappings. They considered several notions of  distortions, which measure the deviation of a map from being distance-preserving (up to a constant). We shall recall one  of these notions, which originates in a paper by Milnor titled  \emph{A problem in cartography} \cite{Milnor}. We recall the precise definitions.

Let $\Omega$ be an open subset of the sphere and let $f \colon \Omega \to \mathbb{E}^2$ be a $C^1$-embedding from $\Omega$ to the Euclidean plane. 
We think of $f$ as being our geographic map.

For a point $x \in \Omega$, the authors in \cite{MO}, following Milnor, defined  two constants, $M$ and $m$ by
 $$M(x):=\sup_{v} |df(v)|/|v|$$ and $$m(x):=\inf_{v}|df(v)|/|v|$$
 where the supremum and the infimum are taken over the set of non-zero tangent vectors $v$ at $x$.
%
%
They defined then the \emph{infinitesimal bi-Lipschitz constant} $\sigma(x)$  of $f$ at $x$ as  $$\sigma(x):=\max\{M(x), m(x)^{-1}\}.$$

Then, they introduced the quantities $$L(f):=\sup_{x \in \Omega} \sigma(x)$$ and $$\ell(f):=\inf_{x \in \Omega} \sigma(x)$$
and they called $$v(f)=\log\left(\frac{L(f)}{\ell(f)}\right)$$ the \emph{metrical distortion} of $f$.

 Using modern computational tools, the authors in \cite{MO}  compared this metrical and the quasiconformal distortions of several maps which they call conical maps, defined on spherical annuli, that is, regions contained between two parallels on the sphere, whose width is in the range of latitudes that Euler considers. The maps include the Delisle--Euler map.  They   
 showed that among the given collection of geographical maps, the Delisle--Euler map is, in a precise sense, the
best one from the metrical  and from the quasiconformal distortions points of view and from others.
 
\section{In guise of a conclusion}

This is not the unique occurrence where a craftsman (here, a geographer, namely, Delisle) has proposed, on the basis of experience or intuition, a heuristic construction or object that satisfies certain optimality conditions, and that it is only many years (sometimes centuries) later that mathematicians proved that, indeed, the objects proposed or used during all this period are indeed optimal. Staying in the field of geography, let us also recall the example of the equidistant azimuthal projection, used for several centuries in the construction of geographical maps,\footnote{In particular this is the map that represents the whole world, on the United Nations' flag.} and whose optimality, from the point of view of metric distortion, among all conformal geographical maps, was proved only in 1969, by Milnor. See the article \cite{Milnor}  and Milnor's comments on this article in his \emph{Collected Works} \cite{Milnor-Collected}.

\end{document}